\documentclass[a4paper,12pt]{article}
\usepackage[utf8]{inputenc}
\usepackage[T2A]{fontenc}  
\usepackage[english]{babel}
\usepackage{hyperref}
\usepackage[centertags]{amsmath}
\usepackage{amsfonts}
\usepackage{amssymb}
\usepackage{amsthm}
\usepackage{graphicx}
\usepackage{xcolor}

\theoremstyle{remark}

\newtheorem*{acknowledgments}{Acknowledgments}
\theoremstyle{definition}
\newtheorem{definition}{Definition}

\newtheorem{example}{Example}
\newtheorem{problem}{Problem}
\usepackage{amssymb}
\usepackage{algorithm}
\usepackage{algorithmic}
\usepackage{sagetex}

\title{On the computation of fundamental functions and Abelian differentials of the third kind}
\author{Yu~Ying*, E.\,A.~Ayryan** ***, 
M.\,D.~Malykh** ****, L.\,A.~Sevastianov** ****}
\date{\small
* Department of Algebra and Geometry, Kaili University,\\
3 Kaiyuan Road, Kaili, China, 556011.\\
\vspace*{0.5 cm}
** Joint Institute for Nuclear Research (Dubna),\\
Joliot-Curie, 6, Dubna, Moscow region,
Russia, 141980.\\
\vspace*{0.5 cm}
*** Dubna State University,\\
19 Universitetskaya St, Dubna, Russia, 141980, Russia\\
\vspace*{0.5 cm}
**** Department of Applied Probability and Informatics, \\
Peoples' Friendship University of Russia (RUDN University),\\
6 Miklukho-Maklaya St, Moscow, Russia, 117198.\\
} 

\begin{document}\large\maketitle

\begin{abstract}
We consider the construction of the fundamental function and Abelian differentials of the third kind on a plane algebraic curve over the field of complex numbers that has no singular points. The algorithm for constructing differentials of the third kind is described in Weierstrass's Lectures. The article discusses its implementation in the Sage computer algebra system. The specificity of this algorithm, as well as the very concept of the differential of the third kind, implies the use of not only rational numbers, but also algebraic ones, even when the equation of the curve has integer coefficients. Sage has a built-in algebraic number field tool that allows implementing Weierstrass's algorithm almost verbatim. The simplest example of an elliptic curve shows that it requires too many resources, going far beyond the capabilities of an office computer. Then the symmetrization of the method is proposed and implemented, which solves the problem and allows significant economy of resources. The algorithm for constructing a differential of the third kind is used to find the value of the fundamental function according to the duality principle. Examples explored in the Sage system are provided.

\end{abstract}

\section*{Introduction}

First successes of computer algebra systems (CAS) are primarily due to the fact that back in the 1960s, it was possible to create algorithms for the symbolic integration of elementary functions, based on Liouville work 
\cite{Moses,Bronstein,Parisse}.  
In the 1980s, the question arose about the development of computer programs that would allow symbolic integration of algebraic functions (Abelian integrals), see 
\cite{Baker}.  In Refs. \cite{Devenport}, \cite{Trager} and \cite{Bronstein-1998} algorithms were developed algorithms for  integration of Abelian integrals  in elementary functions, but the final set of developed algorithms are too complicated to be implemented in CAS  
\cite{malykh-2020-jsc}. 

Of all the known approaches to Abelian integrals, Weierstrass's approach was the most constructive. In Ref. 
\cite{malykh-2020-jsc}, we tried to show that the normal form of representation of Abelian integrals proposed in the lectures gives solutions to a number of classical problems and its implementation in computer algebra systems would be very useful. The key problem on this way, both in the 19th century and now, is the construction of the fundamental function (Hauptfuktion) or, which is also due to the duality principle, the differential of the third kind (Art), the construction algorithm of which is described in the last chapter of Part 1 of the Weierstrass Lectures 
\cite{Weierstrass}, published in 1902 by Hettner and Knoblauch. There are no examples of using the algorithm in the text. 
	
A characteristic feature of Weierstrass' approach is the use of a lot of irrational numbers, the algorithm for determining which is either described in the text, or more or less obvious
\cite{Kochina-1985}. The Sage system 
\cite{Sage} has a built-in implementation QQbar of the  field of algebraic numbers, added by Carl Witty in 2007, so in theory the algorithms from the Lectures can be implemented as written. However, in practice, symbolic expressions containing a ten of numerical coefficients from the field of algebraic numbers QQbar are very difficult to  manipulate. We decided to consider this direct implementation of the algorithms and  these expressions themselves and evaluate the difficulties that arise, see also 
\cite{PCA2018}. 

\section{Abelian differentials of the first kind}

Let polynomial $f$ define an algebraic curve $C$ of the order $r$ on the projective plane $xy$ over the field $\mathbb{C}$. Let for simplicity this curve have no singular points.
	
\begin{definition}
A differential of the form $udx$, $u \in \mathbb{C}(x,y)$, having no singular points on the curve $f$ is called a differential of the first kind.
\end{definition}	

\begin{problem}
\label{p:1}
Given a polynomial $f \in \mathbb{Q}[x,y]$, find a non-constant rational function  $u \in \overline{\mathbb{C}}(x,y)$ such that $udx$ is a differential of the first kind.
\end{problem}
	
The absence of finite singular points makes one seek the solution in the form
	\[
	\frac{E(x,y)dx}{f_y(x,y)}, \quad E \in \mathbb{C}[x,y],
	\]
and the absence of singular points at infinity indicates the fact that the order of the polynomial  $E$ cannot exceed  $r-3$. Since no limitations should be imposed on the coefficients of this polynomial, the set of differentials of the first kind has the dimension
	\[
	p=\frac{(r-1)(r-2)}{2},
	\]
which is called a genus of the curve. For the basis of this space one can take differentials with the coefficients form the field  $\mathbb{Q}$, rather than from its algebraic closure. Therefore, when constructing differentials of the first kind it is possible and necessary to work over the field  $\mathbb{Q}$. 
	
Algorithms for calculating a basis for the space of differentials of the first kind for planar curves, including those having singular points, have been proposed both in classical books and in present-day papers  
\cite{Hoeij}. At present they are implemented in Maple system (AlgCurves, CASA) and partially in Sage. 

\section{Abelian differentials of the third kind}

\begin{definition}
\label{def:iii}
A differential of the form $udx$, $u \in \mathbb{C}(x,y)$ is called a differential of the third kind, if it has two singular points, namely, poles of the first order $(x_1,y_1)$ and $(x_2,y_2)$ with residues  $1$ and $-1$. 
\end{definition}

\begin{problem}
\label{p:2}
Given an indecomposable polynomial $f \in \mathbb{Q}[x,y]$, defining a projective curve $C$, and two points $(x_1,y_1)$ and $(x_2,y_2)$ on this curve, and $x_1,x_2,y_1,y_2 \in \overline{\mathbb{Q}}$.  It is required to construct a non-constant rational function  $u \in \mathbb{C}(x,y)$ such that $udx$ is a differential of the third kind with the poles  $(x_1,y_1)$ and $(x_2,y_2)$.
\end{problem}
	
The addition to the differential of a linear combination of differentials of the first kind does not give rise to new singularities of change of residues, therefore, the solution of Problem \ref{p:2} is defined to a linear combination of $p$ differentials of the first kind. 
	
We describe here briefly the solution to this problem, following 
\cite[ch. 8]{Weierstrass}. The absence of finite singular points with $x\not = x_i$ makes one seek the solution in the form
	\[
	\frac{E(x,y)dx}{(x-x_1)(x_2-x)f_y(x,y)}, \quad E \in \mathbb{C}[x,y],
	\]
and the absence of points at infinity indicates the fact that the order of the polynomial  $E$ cannot exceed $r-1$. Equation
	\[
	f(x_i,y)=0
	\]
beside the root $y=y_i$ has $r-1$ more roots; let us denote them as $y_i', \dots, y_i^{(r-1)}$. If there are no multiple roots among them, then the equations
\begin{equation}
\label{eq:E:1}
E(x_i, y_i^{(j)})=0, \quad i=1,2, \, j=1, \dots, r-1
\end{equation}
ensure the absence of singularities at point, different from  $(x_1,y_1)$ and $(x_2,y_2)$. 

Since there are no singular points on the curve and there are no multiples among the roots $ f (x_1,y) = 0 $, the neighborhood of the point $ (x_1,y_1) $ can be uniformized as 
\[
x=x_1 +t, \quad y=y_1+ c_1t+c_2t^2+\dots
\] 
Substituting these expressions in
\[
\frac{E(x,y)dx}{(x-x_1)(x_2-x)f_y(x,y)},
\]
we have
\[
\frac{E(x_1,y_1)}{(x_2-x_1)f_y(x_1,y_1)} \frac{dt}{t} + \dots
\]
Therefore, the residues at point $ (x_1,y_1) $ will be equal to $1 $ if and only if
\[
E(x_1,y_1)=(x_2-x_1)f_y(x_1,y_1).
\]
Similarly, in the neighborhood $ (x_2,y_2) $ the differential
\[
\frac{E(x,y)dx}{(x-x_1)(x_2-x)f_y(x,y)}
\]
can be described as
\[
-\frac{E(x_2,y_2)}{(x_2-x_1)f_y(x_2,y_2)} \frac{dt}{t} + \dots
\]
Therefore, the residues at point $(x_2,y_2)$  will be equal to $-1$  if and only if
\[
E(x_2,y_2)=(x_2-x_1)f_y(x_2,y_2).
\]
Thus thee conditions for residues at these points give two more equations: 
\begin{equation}
\label{eq:E:2}
	E(x_1,y_1)=(x_2-x_1)f_y(x_1,y_1),
	\quad
	E(x_2,y_2)=(x_2-x_1)f_y(x_2,y_2).
\end{equation}

As a result, the solution to Problem \ref{p:2} reduces to the solution of a system of linear equations \eqref{eq:E:1},\eqref{eq:E:2} with coefficients from \verb|QQbar|, and the main difference of Problem \ref{p:2} from Problem \ref{p:1} is the necessity to extend the number field. 

We  wrote an algorithm for solving the problem \ref{p:2} (Algorithm \ref{alg:iii}), realized it in Sage 
\cite{weierstrass.sage}.

\begin{algorithm}
\caption{Algorithm for finding the third kind integrals in the case of simple roots} 
\label{alg:iii} 
\begin{algorithmic}
\REQUIRE ~~\\ 
   \textbf{Input:} $f\in \mathbb{Q}[x,y]$, points $(x_1,y_1)$, $(x_2,y_2)$ of the curve $f$ over the field $\overline{\mathbb{Q}}$\\
    \textbf{Output: } differentials $udx$ with $p$ arbitrary coefficients, here $u$ is an element of the field of fractions for the ring $\overline{\mathbb{Q}}[x,y]/(f)$. \\
\ENSURE ~~\\ 
   \STATE  \textbf{step1: }  Calculate the lists $R_1$ and $R_2$ of the roots of the equations $f(x_1,y)=0$ and $f(x_2,y)=0$ with respect to $y$. Delete form they the roots $y_1$ and $y_2$. 
   \STATE  \textbf{step2: } Add symbolic variables $c_{ij}$ and define the expression 
   \[
   E=\sum \limits_{i+j\leq r} c_{ij}x^iy^j.
   \]
   \STATE  \textbf{step3: } Calculate the lists  $L$ of the equations 
   \[
   E|_{x=x_i, y \in R_i}=0, \quad i=1,2.
   \]
   Add to it the following  two equations 
   \[
   E(x_1,y_1)=(x_2-x_1)f_y(x_1,y_1), \quad E(x_2,y_2)=(x_2-x_1)f_y(x_2,y_2).
   \]
   \STATE  \textbf{step3bis (optional):} Symmetrization of the obtained equations.
   \STATE  \textbf{step4: } Solve the equations $ L $ with respect to $ c_ {ij} $. Substitute the solution in the expression $ E $ and return
   \[
   \frac{E(x,y)dx}{(x-x_1)(x_2-x)f_y(x,y)}.
   \]
\end{algorithmic}
\end{algorithm}

\begin{example}
\label{ex:iii}
Consider an elliptical curve
\begin{equation}
\label{eq:curve}
x^3-y^3+2xy+x-2y+1=0.
\end{equation}
and construct the differential of the third kind with pols at points with abscissa $x=0$ and $x=1$. 
\begin{sagecommandline}
sage: var("x,y,dx")
sage: f=x^3-y^3+2*x*y+x-2*y+1
sage: x1=0
sage: x2=1
sage: y1=QQbar[y](f.subs(x=x1)).roots(multiplicities=False)[0]
sage: y2=QQbar[y](f.subs(x=x2)).roots(multiplicities=False)[0]
\end{sagecommandline}
First three steps of Algorithm \ref{alg:iii} was realized as function \verb|iii_eqs|, which return six linear equations with six unknowns $c_0, \dots, c_5$:
\begin{sagecommandline}
sage: load("iii.sage")
sage: iii_eqs(f,[x1,y1],[x2,y2])
\end{sagecommandline}
To solve system of equations, Sage uses a standard function \verb|solve|, which does not support the operation with algebraic numbers Therefore, we proceeded to matrices over the field of algebraic numbers and tried to solve the system of linear equations by means of function  \verb|solve_right|.  However, this function did not cope with this system in a reasonable amount of time, returning the warning: increasing stack size to 32000000. 
\end{example}

The above example shows that even in the simplest case, a direct implementation of Algorithm \ref{alg:iii} is very costly. Fortunately, the system of equations \eqref{eq:E:1},\eqref{eq:E:2} consists of two subsystems of the form
\begin{equation}
\label{eq:iii:eqs-old}
E(x_i,y_i^{(j)};c_0, \dots)=b_{i,j}, \quad j=1,2, \dots r,
\end{equation}
where $y_i^{(j)}$ is the set of roots of equation $f(x_i,y)=0$ with respect to $y$. 
Even without specifying the form of the right-hand side, we can obtain from this system $r$ consequences that are symmetric with respect to permutations of the roots of the equation $f(x_i,y)=0$:
\begin{equation}
\label{eq:iii:eqs-new}
\sum \limits_{j=1}^r (y_i^{(j)})^k E(x_i,y_i^{(j)};c_0, \dots)=b_{i,j}, \quad k=0,1, \dots r-1.
\end{equation}
The coefficients at unknowns in the new system are symmetric functions of the roots, therefore, they will be rational numbers if such is the considered value of $ x_i $. Therefore, to solve the new system, it is sufficient to invert the matrix with rational coefficients. System \eqref{eq:iii:eqs-new} is equivalent to original system \eqref{eq:iii:eqs-old}, since the new system is obtained from the old one by multiplying by the Vandermonde matrix
\[
\begin{pmatrix}
1 & \cdots & 1\\
\vdots & \ddots & \vdots \\
(y_i^{(1)})^{r-1} & \cdots & (y_i^{(r)})^{r-1}
\end{pmatrix},
\]
the determinant of which in the case of simple roots considered is nonzero. The described transition to symmetrized subsystems will be referred to as SLAE symmetrization and always carried out before using the Gauss method (step 3bis in the algorithm \ref{alg:iii}).
\begin{example}
\label{ex:iii:bis}
Let us return to the example \ref{ex:iii}. After implementing the first three steps and the symmetrization in the form of function \verb|iii_eqs_sym|, we obtain six linear equations with six unknowns $c_0, \dots, c_5$:
\begin{sagecommandline}
sage: iii_eqs_sym(f,[x1,y1],[x2,y2])
\end{sagecommandline}
Now the drawbacks of in the realization of the field of algebraic numbers are obvious: one of coefficients in the second equation is not identified as rational. Note, that when operating with algebraic numbers, there is no rounding error and you can verify its rationality using standard tools:
\begin{sagecommandline}
sage: eqs=iii_eqs_sym(f,[x1,y1],[x2,y2])
sage: eqs[1]
sage: eqs[1].coefficient(c2)
\end{sagecommandline}
The symmetrized SLAE is solved without noticeable expenditure of time. For convenience we constructed our own user function \verb|lsolve| which  reduces the equations to the matrix form and  solves them by means of function \verb|solve_right|:
\begin{sagecommandline}
sage: lsolve(eqs,[c0,c1,c2,c3,c4,c5])
\end{sagecommandline}
This function, like the function \verb|solve_right| itself, returns a partial solution. The entire algorithm \ref{alg:iii} is implemented as function \verb|iii| that returns the differential of the third kind, defined up to a linear combination of differentials of the first kind:
\begin{sagecommandline}
sage: iii(f,[x1,y1],[x2,y2])
\end{sagecommandline}
In the example considered, a visually graspable expression is obtained.
\end{example}
Thus, such symmetrization is quite enough for efficient implementation of the method for constructing a differential of the third kind, proposed in Weierstrass's Lectures.
\section{Fundamental function}
The theory of algebraic curves is based on a seemingly very simple problem.
\begin{problem}
\label{p:r}
Given an indecomposable polynomial $f\in\mathbb{C}[x, y] $ defining a projective curve $ C $, and $ s $ points $ (x_i, y_i) $ on this curve. It is required to construct a non-constant rational function $ g \in \mathbb {C} (x, y) $, having poles only at given points and only of the first order.
\end{problem}
If the points are not chosen in a special way, then this problem is solvable only if the number $ r $ is greater than some boundary $ p $, remarkably equal to the genus of the curve.
\begin{definition}
The points $ (x_i, y_i) $ of the curve $ C $ will be called points in general position if the problem \ref{p:r} has no solution.
\end{definition}
It is easy to show that solving the problem \ref{p:r} for any $ s $ reduces to solving the problem for $ s = p + 1 $, while the solutions of the problem \ref{p:r} for points $ (a_1, b_1), \dots, (a_p, b_p) $ and $ (x ', y') $ are defined up to two constants, additive and multiplicative. If we fix the multiplicative constant for the selection of the residue at the point $ (x ', y') $, then the solution will be determined up to an additive constant. This solution is referred to as fundamental function in Weierstrass's lectures.
\begin{definition}
\label{def:haupt}
We fix $ p $ points $ (a_i, b_i) $ in general position on the curve $ C $. The fundamental function is the solution of the problem \eqref{p:r} for the set of points obtained by adding one more point to these $ p $ points, which will be denoted below as $ (x ', y') $. For normalization, it is assumed that the residue at this point is equal to $ -1 $.
\end{definition}
The definition \ref{def:haupt} is convenient for proving the existence of a fundamental function
\cite[ch. 2]{Weierstrass}.
The question of its direct use for finding the fundamental function on a given curve has apparently never even been raised. In the Lectures, the key to both the use of Abelian integrals in the theory and to the construction of the fundamental function is the connection between the fundamental function and the integral of the third kind.
Let $ udx $ be an Abelian differential of the third kind, i.e., a differential (def.~\ref{def:iii}) having two singular points, namely, first-order poles $ (x_1, y_1) $ and $ (x_2, y_2 ) $ with residues $ 1 $ and $ -1 $, and let $ h $ be the fundamental function (def.\ref{def:haupt}). Then the product $ hudx $ has poles of the first order at $ p + 3 $ points and the residue theorem yields
\[
-h|_{(x_1,y_1)}+h|_{(x_2,y_2)} + u|_{(x',y')} +\sum \limits_{i=1}^p c_i u|_{(a_i,b_i)}=0,
\]
where $c_i$ is the residue of the fundamental function at the point $(a_i,b_i)$. Since the fundamental function is defined up to an additive constant, it is possible to assume that $h$ is zero at the point $(x_2,y_2)$, then
\[
h|_{(x_1,y_1)} = u|_{(x',y')} +\sum  \limits_{i=1}^p c_i u|_{(a_i,b_i)}.
\]
Since a differential of the third kind is defined up to $p$ constants, it is always possible to ensure that $u|_{(a_i,b_i)}=0$. But then
\begin{equation}
\label{eq:hu}
h|_{(x_1,y_1)} = u|_{(x',y')},
\end{equation}
i.e., the fundamental function multiplied by $dx'$, as a function of the point $(x',y')$ is a differential of the third kind. In Weierstrass’s lectures a direct proof is given to this fact, from which both the equality of the genus and the dimension of the space of differentials of the first kind, and the existence of differentials of the third kind immediately follows.
\begin{algorithm}
\caption{ Algorithm for calculating the value of the fundametal function at a given point $(x_1,y_1)$}
\label{alg:hu}
\begin{algorithmic}
\REQUIRE ~~\\ 
   \textbf{Input:} $f\in \mathbb{Q}[x,y]$, points $(x_1,y_1)$, $(x_2,y_2)$,$(x',y'),(a_1,b_1), \dots, (a_p,b_p)$ of curve $f$ over field $\overline{\mathbb{Q}}$\\
    \textbf{Output: } value of the fundamental function $h$ at point $(x_1,y_1)$, having simple poles at points $(x',y'),(a_1,b_1), \dots, (a_p,b_p)$, residue $-1$ at point $(x',y')$ and zero at point $(x_2,y_2)$. \\
\ENSURE ~~\\ 
   \STATE  \textbf{step1: } Using points $(x_1,y_1)$, $(x_2,y_2)$ construct differential $udx$ of the third kind, defined up to $p$ constants:
   \[
   u + \sum \limits_{i=1}^p c_i u_i,
   \]
   where $u_i$ are integrals of the first kind.
   \STATE  \textbf{step2: } Determine these constants from $p$ equalities
   \[
   u|_{(a_i,b_i)}  + \sum \limits_{j=1}^p u_j|_{(a_i,b_i)} c_j  =0, \quad i=1,2,\dots, p.
   \]
   \STATE  \textbf{step3: } Return
   \[
   \left. u + \sum \limits_{i=1}^p c_i u_i \right|_{(x',y')}
   \]
\end{algorithmic}
\end{algorithm}
For us now, it is more important that the formula \eqref{eq:hu} allows calculating the value of the fundamental function at almost any point according to the algorithm \ref{alg:hu}. It should be noted that the first step can be performed according to the algorithm \ref{alg:iii} only if the equations $ f (x_0, y) = 0 $ and $ f (x_1, y) = 0 $ do not have multiple roots, i.e., we can calculate the values of the fundamental function not at all points, but at almost all points.
\begin{example}
Consider again the elliptic curve
\[
x^3-y^3+2xy+x-2y+1=0
\]
of genus 1, at which we take the point $ (1, y_2) $ as a zero of the fundamental function, and the point $ (2, b_1) $ as an additional pole $ (a_1, b_1) $. Find the value of the fundamental function at point $ (0, y_1) $. We define the ordinates of the points so that they fall on the curve
\begin{sagecommandline}
sage: f=x^3-y^3+2*x*y+x-2*y+1
sage: x1=0
sage: x2=1
sage: a1=2
sage: xx=3
sage: y1=QQbar[y](f.subs(x=x1)).roots(multiplicities=False)[0]
sage: y2=QQbar[y](f.subs(x=x2)).roots(multiplicities=False)[0]
sage: b1=QQbar[y](f.subs(x=a1)).roots(multiplicities=False)[0]
sage: yy=QQbar[y](f.subs(x=xx)).roots(multiplicities=False)[0]
\end{sagecommandline}
Function \verb|haupt_fuction_eval| returns the value of the fundamental function at point $(x_1,y_1)$:
\begin{sagecommandline}
sage: haupt_fuction_eval(f,[x1,y1],[x2,y2],[xx,yy],[[a1,b1]])
\end{sagecommandline}
Unfortunately, a very unexpected difficulty arises here: the value of the fundamental function is an algebraic number, the minimum polynomial for which cannot be calculated in a reasonable time (the warning is
“increasing stack size to 64000000”).
\end{example}
In the process of designing the function \verb|haupt_fuction_eval| we faced a number of problems that led to the same warning and offered no possibility to finish the calculations. We managed to avoid them by dividing a number of symbolic expressions into parts, each calculated separately. We had to apply the method \verb|subs| to individual elements rather than to lists.

\section{Conclusion}
Summarizing the above, it can be argued that the implication of the algebraic number field in Sage really allows, at least in simple examples, to implement Weierstrass's algorithms almost as described in his Lectures. In doing so, however, it is important to carry out symmetrization wherever possible. Otherwise, it will not be possible to get a solution within a reasonable time.

The next step in implementing algorithms, proposed in Weierstrass's Lectures, is the construction of the fundamental function. For this purpose, it is sufficient to construct a differential of the third kind with a movable pole. To execute symmetrization in this case, too, we intend to use a perfect tool  --- the package Symmetric Functions for Sage, which allows expressing a symmetric function from a ring $K[x_1, \dots, x_n]$  as a linear combination of elementary symmetric functions.  

\begin{acknowledgments}
The computations presented in the paper, were performed in Sage (\url{www.sagemath.org}). The contribution of  E.\,A.~Ayryan (Investigation), M.\,D.~Malykh (Investigation, algorithms development), and L.\,A.~Sevastianov (Conceptualization, writing) is supported by the Russian Science Foundation (grant no. 20-11-20257).
\end{acknowledgments}

\bibliographystyle{ugost2008}
\bibliography{malykh}

\end{document}